\documentclass[12pt,reqno]{amsart}

\newtheorem{theorem}{Theorem}[section]
\newtheorem{lemma}[theorem]{Lemma}

\newtheorem{proposition}[theorem]{Proposition}

\newtheorem{example}[theorem]{Example}

\begin{document}

\title[Diffeomorphisms with Holes]
{Topological and Symbolic Dynamics for Axiom A Diffeomorphisms with Holes}
\author{Stefan Bundfuss}
\address{FB Mathematik, Technische Universit\"at Berlin,
Germany} 
\email{bundfuss@math.tu-berlin.de}
\urladdr{http://www.math.tu-berlin.de/{\lower.7ex\hbox{\~{}}}bundfuss}
\author{Tyll Kr\"uger}
\address{FB Mathematik, Technische Universit\"at Berlin and Fakult\"at
f\"ur Physik, Universit\"at Bielefeld,
Germany} 
\email{tkrueger@math.tu-berlin.de}
\author{Serge Troubetzkoy}
\address{Centre de Physique Theorique,
F\'ed\'eration de Recherche des Unit\'es Math\'ematiques de Marseille,
Institut de Math\'ematique Luminy, and
Universit\'e de la M\'editerran\'ee, Marseille, France}
\email{serge@cpt.univ-mrs.fr}
\email{troubetz@iml.univ-mrs.fr}
\urladdr{http://iml.univ-mrs.fr/{\lower.7ex\hbox{\~{}}}troubetz}
\date{\today}
\thanks{This paper was written during reciprocal visits which where 
supported by SFB 288, PRODYN and Universit\"at Bielefeld}



\begin{abstract}
We consider an Axiom A diffeomorphism and the invariant
set of orbits which never falls into a fixed hole.  We 
study various aspects of the complexity of the symbolic
representation of $\Omega$.  Our main result are that each
topologically transitive component of $\Omega$ is coded
and that typically $\Omega$ is of finite type.
\end{abstract}

\maketitle

\renewcommand{\theequation}{\arabic{section}.\arabic{equation}}
\def\w{\omega}
\def\hw{\hat{\w}}
\def\g{\gamma}
\def\a{\alpha}
\def\ba{\begin{eqnarray}}
\def\ea{\end{eqnarray}}
\def\Se{\S_{\text{even}}}
\def\be{\begin{equation}}
\def\ee{\end{equation}}
\def\nn{\nonumber}
\def\O{\Omega}
\def\T{{\bf T}}
\def\S{\Sigma}
\def\hS{\hat{\Sigma}}
\def\iH{H}
\def\bH{\partial H}
\def\bB{\partial B}
\def\d{\text{dist}}
\def\P{{\mathcal P}}
\def\Q{{\mathcal Q}}
\def\e{\epsilon}
\def\ux{\underline{x}}
\def\H{\hat{H}}
\def\hx{\hat{x}}
\def\hy{\hat{y}}
\def\hc{\hat{c}}
\def\hz{\hat{z}}
\def\e{\varepsilon}

\section{Introduction}\label{secI}
Let $f$ be an Axiom A diffeomorphism of a compact $s$--dimensional 
Riemannian manifold and $\Lambda \subset M$ the nonwandering set of $f$.
We fix a Markov partition, this yields a 
natural representation of $\Lambda$ as a subshift
of finite type $\hS$. We cut out an open hole $H$ out of $M$ and 
consider the invariant set $\O = \O_H \subset \Lambda$ of nonwandering 
points whose 
orbit (forward and backwards) never falls in the hole.  
The set $\O$ corresponds to a subshift $\S = \S_H
\subset \hS.$  We call such a subshift $\S$ an {\em exclusion subshift}.
We are interested in several questions about the topological structure
of exclusion subshifts.  More precisely for which holes are exclusion
shifts  subshifts
of finite type (SFT), sofic shifts, or coded systems?
As we vary the hole what is the typical type of an exclusion subshift?

Our first result is that an exclusion subshift
has at most countably many topologically transitive components
and exclusion subshifts are always coded systems on
each topological transitive component. Under additional assumptions 
we can prove the finiteness of the number of transitive components.
Next we give a criterion for when the hole leads to a SFT.
We then apply this criterion to several different classes of holes
(rectangles, polyhedra, holes with continuous boundary) with corresponding
topologies/measures to show
that the ``typical'' hole leads to a SFT.
Finally we show that every $\beta$--shift is an exclusion shift, thus
in particular there are exclusion shifts which are sofic and ones which
are not sofic.

Dynamical systems with holes have been studied extensively by phy\-sicists.
Cvitanovich and his coworkers have investigated how to 
characterize the set $\S$ in various settings
(see the survey \cite{Cv}). They have
introduced the notion of a pruning front which corresponds to
$\bH$ in our setting.  Some aspects of Cvitanovich's work
have been carried out by Carvalho
in a mathematical  framework for Smale's horseshoe \cite{Ca}. 
Another kind
of question about dynamical systems with holes have also been
extensively studied by physicists, namely 
construction of a physical semi--invariant measure and the understanding
of the speed of mass disappearance
into the holes (the escape rate formula) \cite{GD}.   
A series of mathematical works have
confirmed the expectations of the physicists in many settings
\cite{C1,C2,CM1,CM2,CMT1,CMT2,CMS,LM,PY,R}.
B\"aker and Kr\"uger have initiated the study of the topological structure
of exclusion subshifts in the one dimensional setting \cite{BK}.  

\section{Definitions}\label{def}

Let $M$ be a compact $s$--dimensional ($s \ge 2$) 
Riemannian manifold and $f: M \to M$ an Axiom A diffeomorphism.
Fix a proper\footnote{A Markov partition is called proper if each element
of the partition is the closure of its interior.}, generating Markov 
partition\footnote{More generally, we can consider a map which is a
local diffeomorphisms except on a singularity set which consists of
a finite union of codimension one manifolds such that the map
admits a finite proper generating Markov partition.}. 
Let $(\hS,\sigma)$ be the resulting SFT and $\pi: \hS \to M$
the projection map.  We cut out an open\footnote{Occasionally we will
discuss cases when the hole is not open.} hole out of $M$ 
whose boundary
$\bH = \bar{H}\backslash H$ consists of a finite union of topological
$(s-1)$--dimensional spheres and $H = int(\overline{H})$. 
Consider the invariant set 
$\O^* = \O^*_H$ of points whose 
orbit (forward and backwards) never falls in the hole
and the invariant set $\O$ of nonwandering points
for $f|_{\O^*}$.  
Let $\S^* = \S^*_H = \pi^{-1} \O^*$ and $\S = \pi^{-1} \O$ 
with the following convention\footnote{This convention is not needed in the
case that the coding is always unique, for example for horseshoes}:  
if $x \in \O$ is on the boundary of
the hole and the preimage of $x$ is not unique, then we only
consider those preimages of $x$ to be in $\S$ which can be approximated
from outside the closure of $H$, i.e.~those $s \in \pi^{-1}(x)$ 
such that $\exists x_j \not \in \overline{H}$
such that $x_j \to x$ and $\exists s_j \in \pi^{-1}(x_j)$ with $s_j \to s$.
We call such a subshift $\S$ an {\em exclusion subshift}.
Remark:  our standard examples will be two dimensional,
the $n$--branch horseshoe map and
the $n$--branch Bakers map.

If $M$ is two dimensional, 
call a hole a {\em rectangle like hole} if $\bH$ consists of
a finite number of curves, each of which is
a finite length piece of stable or unstable manifold of $f$.  We call the
corresponding exclusion shift a {\em rectangle exclusion shift} RES.

We will also consider the one dimensional situation,
where $M$ is an interval or a circle. We will consider continuous maps
of the circle or piecewise continuous maps of the interval which, when
considered as a map of the circle are continuous (for example the doubling
map $f(x) = 2x $ mod 1). Additional we require the existence of a finite 
generating Markov
partition $\P$ whose elements consist of intervals.  
Since the map is not invertible we only require that the forward orbit never
falls into a hole.  If the hole consists
of a finite union of intervals we call the corresponding exclusion
shift an {\em interval exclusion shift} IES.  Our standard example
in this framework is the doubling map.

\section{SFTs and sofic systems}

Consider the alphabet $\{1,2,\dots,n\}$ and a finite collection
of forbidden words on length $m$ ($m$ fixed, usually taken to be 2).
The subset of all sequences in $\{1,\dots,n\}^{\mathbb{Z}}$ 
(resp.~$\{1,\dots,n\}^{\mathbb{N}}$ )
where the forbidden word never appears is called
a {\em subshift of finite type} (SFT).

\begin{theorem} Every SFT is an exclusion subshift.
\label{represent}
\end{theorem}

\begin{proof}
Consider a SFT $\S$.  Let $n$ be the cardinality of the alphabet of $\S$
and $m$ the length of the forbidden blocks.
Consider a horseshoe map $f$ with $n$--branches. 
We consider the standard Markov partition for $f$.
Define $H'$ to be the (finite) union of Markov rectangles (for $f^m$)
which correspond the forbidden $m$--blocks which define $\S$.
Markov rectangles are by definition closed, however 
since a horseshoe is a Cantor set we can find an open hole $H$ containing
$H'$ such that the intersection of $H$ with the horseshoe is exactly $H'$. 
This yields $\S$ as an exclusion subshift. 
\end{proof}

\noindent
Remark: the same construction works for the $n$-fold Bakers map, except that
we define the hole $H$ to be the interior of $H'$ and use the coding
convention from section \ref{def}.

Let $\Se^+$ be the one sided
even shift, that is the set of all 0--1 half infinite
sequences with the constraint that the number of consecutive  ones which
occur in between two zeros is always even.  
  
\begin{example}
\label{example1}
Let $M:=[0,1]$, $f(x):= 2x \mod 1$, $\P=\{ [0,\frac12], [\frac12, 1]\}$.
The even shift $\Se^+$ is not an interval exclusion system with respect to
$(f,\P).$
\end{example}

\begin{proof}
Suppose we have a hole which is a union of intervals and yields
$\Se^+$ as an exclusion system.
Consider the point $x_n=01^{2n+1}0^{\infty}$. 
The whole orbit of $x_n$, except $x_n$ itself is in $\S^+$. 
Thus  $x_n$ must lie in the hole.

Let $y_n=01^{2n}0^{\infty}$. Clearly $y_n \in \S^+$ and 
$x_{n-1} < y_n < x_n$. 
Thus the hole must consist of an infinite number of intervals:
$\Se^+$ can not be an IES.
\end{proof}

Remark:  the exact same proof shows that $\Se^+$ is not an
exclusion system with respect to $(f,\P)$ for any hole
which is a finite union of intervals which are not necessarily
open.

Let $\Se$ be the two sided
even shift, that is the set of all 0--1 bi--infinite
sequences with the constraint that the number of consecutive  ones which
occur in between two zeros is always an even number.  
\begin{example}
\label{example2}
Let
$M = [0,1]^2$,  $f: M \to M$ be the Bakers map:
$$f(x,y) =\left\{\aligned  (2x \mod 1, y/2)\quad&\text{if}\quad x \le 1/2 \\ (2x \mod 1,(y+1)/2)\quad&\text{otherwise}.\endaligned\right. $$
Let $\P := \{P_0,P_1\}$ be the standard Markov partition, 
i.e.~$P_0 := \{(x,y) : x \le 0.5\}$ and $P_1 := \{(x,y) : x \ge 0.5\}$
Then the even shift $\Se$ is not a RES with respect to $(f,\P)$.
\end{example}

\begin{proof}
Let $x_n=0^{\infty}.1^{2n+1}0^{\infty}$ (here the decimal point marks
the position between the $-1$st and $0$th elements of the sequence).
The point $x_n$ is not in the even shift,
it must fall into the hole under some iteration of $f$. 

We treat several cases, first of all suppose that $x_n$ falls into the hole
at the boundary of $M$. Then at the instant that $x_n$ falls into the hole
all the $1$'s are to the right (or all are to the left) of the decimal point
(i.e. $f^jx_n \in \partial M$ with $j \le 0$ or $j \ge 2n+1$).
Note that the intersection of the rectangular holes with the boundary 
consists of a finite union of intervals.  Thus we can apply the
argumentation of the previous
example to conclude that it is impossible to have an infinite number of
the $x_n$ fall into the hole when they are on the boundary of $M$.

In other words all but finitely many $x_n$ which fall into the hole
away from the boundary of $M$.  For such an $x_n$ consider the code
$a$ at the instance of falling into the hole.  It has the form
$0^{\infty}1^p.1^{q+1}0^{\infty}$ where $p+q=2n$.
 
Consider the sequence $(u_i)_{i\in {\bf Z}}$ with 
$u_i=0^{\infty}1^{2i}01^{p-1}.1^{q+1}0^{\infty}$ 
and the sequence $(v_i)_{i\in {\bf Z}}$ with
$v_i=0^{\infty}1^p.1^q01^{2i}0^{\infty}$.
These two sequences get arbitrary near from the left
and the bottom to $a$ and all their elements are in the even shift. 
Thus $a$ must be a corner of the hole. 
This contradicts to the assumption that the number of corners of the hole 
is finite.
\end{proof}

\noindent
Remark: the same proof shows that $\Se$ is not a RES with respect to the
$n$--fold Bakers map, it also shows that $\Se$ is not a RES 
for any hole which is a finite union of closed rectangular
like holes. An even simpler proof can be given if
$f$ is a standard horseshoe map with a standard
Markov partition.

A closed shift invariant subset of $\{1,2,\dots,n\}^{\mathbb{Z}}$ 
(resp.~$\{1,2,\dots,n\}^{\mathbb{N}}$) 
is called 
{\em sofic} if it a factor of a SFT.  The even shift is an example of a sofic
system.

\begin{theorem}
\begin{enumerate}
\item{In dimension one not every sofic system is an interval exclusion 
subshift.}
\item{In dimension two not every sofic system is a rectangular hole 
exclusion subshift.} 
\end{enumerate}
\end{theorem}

\begin{proof}
In the two examples we showed that the even
shift (one sided even shift) can not be represented as a RES (IES) 
in our standard model: the Bakers map (the doubling map).

To prove part 2 we will show that if the even shift is
representable as a RES for some Axiom A diffeomorphism then it
is in fact an RES for the Bakers map. The proof of part one is similar.

Let $(f,\P)$ be the Bakers map with the standard Markov partition.
We suppose by way of contradiction that $\Se$ is a RES for 
some Axiom A diffeomorphism $g$ with respect to the Markov partition $\Q$.  
Consider the SFT $\S$ defined by $(g,\Q)$
and denote the projection $\pi: \S \to \Lambda_g$.
By assumption there is a rectangle
like hole $H$, that is the set $\O_g$ of points which never fall into $H$
satisfies $\pi: \Se \to \O_g$.

By the proof of theorem \ref{represent} we can find an rectangle like hole
$H'$ which yields $\S$ as a RES for
$(f,\P)$.  Let $\Lambda'_f$ be the invariant set of points which never fall 
into $H'$ and let $\pi': \S \to \Lambda'_f$ be the projection.  

Consider $H_{\text{symb}} := \{ s \in \S: \pi(s) \in H\}$ and
$H'' :=\pi'(H_{\text{symb}})$.
Since the maps $\pi$ and $\pi'$ preserve local product structure the set $H''$
is a rectangle like hole.  If $H''$ is open then $\Se$ is a RES in
the Bakers map, a contradiction.
If $H''$ is not open then the interior of $H''$ is also a rectangular
like hole, which
using the convention defined in section \ref{def} shows that $\Se$ is an
RES  for the Bakers map, a contradiction. 
\end{proof}

We do not know if every sofic system, or even if every coded system is
an exclusion subshift.  However, if we drop the requirement that the
boundary consists of spheres then
we get that {\bf all} subshifts are {\em weak exclusion
subshifts}, that is they are the complement of
an open hole without any further assumptions on the boundary. 
To see this simply take the 
compliment of the invariant set which gives rise to the subshift as
the hole.

\begin{theorem}
Every $\beta$--shift is a RES.
\end{theorem}

\begin{proof}
We begin by describing a one dimensional construction from \cite{BK}.
Consider the dynamical system $h_{\beta}: [0,1) \to [0,1)$ given by
$h_{\beta}x = \beta x \hbox{ mod } 1.$  
Suppose $\beta \in (1,2]$ and consider the 
``partition'' $I_0 := [0,\beta^{-1}]$ and $I_1 := [\beta^{-1},1]$.
Given $\omega,\alpha \in \{0,1\}^{\bf N}$ we say $\omega < \alpha$
if there exists a positive integer $N$ such that 
$\omega_i = \alpha_i$ for $i=0,1,\dots, N-1$ and $\omega_N=0 < \alpha_N=1$. 
It is well known that the code 
$\hat{\beta}$ of the orbit
of $\beta$ satisfies $\sigma^n \hat{\beta} < \hat{\beta}$ for all $n>0.$
The set of all sequences which are the code of some orbit is called the
one sided $\beta$--shift, it is characterized by \cite{B} 
$$\{x \in \{0,1\}^{\bf N}: \sigma^n x < \hat{\beta} \hbox{ for all } n > 0\}.$$
Let
$g : [0,1) \to [0,1)$  be defined by $g(x) = 2 x$ mod $1.$ Inspired by the
above a number $\beta$ is called a $\beta$--number 
if $g^n \beta < \beta$ for all
$n > 0.$  If $\beta$ is a $\beta$--number then the set $X_{\beta} :=
\{x \in [0,1):
g^n x < \beta \ \forall n > 0\}$ is conjugate to 
the one sided $\beta$--shift. If
$H \subset [0,1)$ is the interval $(\beta,1)$ then
the set of points whose forward $g$ orbit never falls in $H$ is exactly 
$X_{\beta}.$

Now let $f$ be the Bakers map.
Fix $\beta$ a $\beta$ number and let 
$H \subset M$ be the rectangular hole $\{(x,y): \ x > \beta, \ 0\le y \le 1\}.$
Let $(x_{-i},y_{-i}) = f^{-i}(x,y).$  Since the hole
stretches from the bottom to the top of $M$ it is easy to see
that $(x,y) \in \O$ if and only if $g^n(x_{-i}) < \beta$ for all
$n \ge 0$ and all $i \ge 0.$  Thus the exclusion shift in this
example  is exactly the natural extension of the one sided $\beta$--shift.
We remark that the two sided $\beta$--shift is of finite type, sofic, or not
sofic  if and only if the one sided one has the same property.  
It is well known that the $\beta$--shift is
of finite type if the binary expansion of $\beta$ is finite, it is sofic
if it eventually periodic and otherwise it is not sofic \cite{B}.
\end{proof}
 
Remark: for $\beta \in (n-1,n]$ we must use the map $g = nx \mod 1$.

\section{Transitive components of $\S$ are coded}

A homeomorphism $T$ of a compact metric space $X$ is called {\em topologically
transitive} if there exists a dense orbit, or equivalently if every proper
closed $T$-invariant subset is nowhere dense.
A closed invariant set $\hat{X}$ (i.e.~$T^{-1}A = TA = A$)
is called a {\em topological 
transitive component}
if  $\hat{X}$ contains a dense orbit and no closed invariant 
set $\hat{X}' \supset \hat{X}$ contains a dense orbit.

We remark that the topologically transitive components of a SFT are
exactly it's irreducible components since, by convention (see section 2) 
we restrict the map to
it's nonwandering set (see for example \cite{LiM}).  In particular
this implies that the topologically transitive components of a SFT
are disjoint.
Furthermore periodic points are dense in a topologically transitive
component of a SFT or in a topologically transitive component of
an Axiom A diffeomorphism.

We call a topologically transitive component $\hat{X}$ {\em boundary
supported} (BSC) if periodic points whose orbit avoids $\overline{H}$
are not dense in $\hat{X}$ otherwise we call $\hat{X}$ {\em non boundary
supported}.

A subshift $\S$ is called a {\em coded} system if it can
be represented by an irreducible countable labeled graph \cite{BH1}.
Equivalently, $\S$ is called coded if $\S$ contains an increasing sequence
of irreducible subshifts of finite type (SFTs) 
whose union is dense in $\S$ \cite{BH2}.

\begin{theorem}\label{toptrans}

Every non boundary supported topologically transitive component is coded. There
are at most countably many such topologically transitive components. 
\label{prev}
\end{theorem}

\begin{proof} Let $\S$ be an exclusion subshift, i.e.~$f$ is an
Axiom A map,  $\P$ a fixed proper, generating 
Markov partition, and $H$ an open hole with compact boundary.   
Then $\P^{(n)} := \vee_{i=-n}^n f^i \P$ 
is also a proper, generating
Markov partition.  Let $H^{(n)} := \text{int}(\cup_{\{ P \in \P^{(n)}: \ 
P \cap \overline{H} \ne \emptyset\}} P).$  Clearly  $H^{(1)} \supset H^{(2)} \supset
\cdots$ and $\overline{H} \subset \cap_n H^{(n)}$.
Furthermore, since $\P$ is proper we have 
$\text{diam}(\P^{(n)}) \to 0$ as $n \to 
\infty$ thus $\cap_n H^{(n)} = \overline{H}.$  
The ``exclusion'' system 
$\S^{(n)} := \S_{H^{(n)}}$ is a SFT.\footnote{We write quotes since the boundary
of $H^{(n)}$ may not be the union of codimension 1 spheres.} 
Clearly $\S^{(1)} \subset \S^{(2)} 
\subset \cdots$
and 
\begin{equation}
\tilde{\S} := \overline{\bigcup_n \S^{(n)}} \subset \S.
\label{rrr}
\end{equation}

Clearly any topologically transitive component which is a subset of
$\overline{\bigcup_n \S^{(n)}}$ is non boundary supported since periodic points
are dense in SFTs.  On the other hand if $x \in \S \setminus 
\tilde{\S}$ then $x$ is not an accumulation point of
periodic points $x_n$ which avoid $\overline{H}$ since any such periodic
point necessarily belongs to $\S^{(n)}$ for some $n$. 
Thus we have shown that every non boundary supported
topologically transitive component
is contained in $\tilde{\S}$.

If each $\S^{(n)}$ was topologically
transitive then $\tilde{\S}$ would be coded. 
Let  $A^{(n)}_i$ be the topologically transitive components of $\S^{(n)}$.
The $A^{(n)}_i$ form a filtration in the sense that
for each $A^{(n)}_i$ there exists $A^{(n+1)}_j$ such that $A^{(n)}_i
\subset A^{(n+1)}_j$. In other words the transitive components $A^{(n)}_i$ 
form an at most countable union of directed trees with each nodes out degree 
is exactly one. The equality in equation $\eqref{rrr}$ implies that each
topologically transitive component of $\sigma|_{\Sigma}$ contains 
a set $\cup A^{(n)}_i$ where the union is taken over a path in one of the
trees (we call this a path limit).
Noticing that such a path is uniquely defined by the root of the
tree since the out degree is always one, implies that there are at most
countably many such paths and thus $\sigma|_{\tilde{\S}}$ has at most countably
many topologically transitive components.
This finished the proof of the countability of the claim of the
theorem.

Suppose $C \subset \tilde{\S}$ is a topologically transitive component of 
$\sigma|_{\tilde{\S}}$.  
To see that $C$ is coded we will define a new filtration.
Since $\S^{(n)}$ is a subshift of finite type it has finitely many
topologically transitive components $A_i^{(n)}$ which are pairwise disjoint.
Consider those components $A_i^{(n)} : i =1,\dots,k_n$ which
are strictly contained in $C$.  
We can assume that the $A_i^{(n)}$
are so ordered that $A_1^{(n)} \subset A_1^{(n+1)}$ for all $n.$
We only need to show that
\begin{equation}\label{=}
\overline{\bigcup_n A_1^{(n)}} = C.
\end{equation}
The rest of the proof is devoted
to establishing this equality.

If equation \eqref{=} is not true there is a $n_0$ such that for all $n \ge n_0$
we can find another $A_i^{(n)} \subset C$ which we denote without loss of
generality 
$A_2^{(n)}$ such that $A_2^{(n)} \subset A_2^{(n+1)}$ 
for all $n \ge n_0$ but ${\cup_n A_2^{(n)}} \cap 
{\cup_n A_1^{(n)}} = \emptyset$. In the terminology introduced above this
means we can find two disjoint paths in the trees whose path limits are both
contained in $C$.

Fix $n \ge n_0.$ For $i=1,2$ consider
a finite word $w_i \in A_i^{(n)}$ where  each symbol and each transition which
characterize $A_i^{(n)}$ appear in $w_i.$ Since $C$ is 
topologically transitive there is a point
$\ux \in C$ where the words $w_1$ and $w_2$ both appear in $\ux.$ Thus we
can find $l$ (which we assume positive without loss of generality) so that
$\ux \in \sigma^l w_1 \cap w_2.$  

Consider the point $x \in M$ with symbolic coding $\ux.$  
We can assume without loss of generality
that $f^i x \not \in \partial H$ for $i=0,\dots,l.$  From our 
assumptions for any sufficiently 
large $N$ and for $i=0,\dots,l$ the rectangles  
$P_i \in \P^{(N)}$  defined by $f^ix \in P_i$ are disjoint
from $H$.
We can form a new SFT,
$\hat{\Sigma}$, with alphabet
the $P_i$ and all the transitions made by the orbit segment $f^ix: i=0,\dots,
l$ allowed. 
This is clearly a topologically transitive SFT.  It also contains $A_i^{(n)}$
for $i=1,2$ since any legal transition in these sets is a legal transition
in $\hat{\Sigma}$.  But, by the construction of $\S^{(m)}$ we have
$\hat{\Sigma} \subset \S^{(m)}$ for sufficiently large $m$ and thus
$\hat{\Sigma} \subset A_1^{(m)}$ and $A_2^{(m)}$ a contradiction. 
\end{proof}

Remark: we actually have a stronger property than coded since our subshifts 
are well approximable from outside as well as inside.

\subsection{Finiteness of topological transitive components in dimension 1}

Suppose that $M$ is one dimensional, $(f,\P)$ is as described in section
\ref{def} and $H = \cup_{i=1}^p H_i$ where 
the $H_i$ are disjoint open 
intervals. The definitions of $\O^*,\O,\S^*$ and $\S$ are similar to
the corresponding definitions
in the invertible case with the only difference being that we
only require that the forward orbit does not fall into $H$.
Let $\O^c := I \backslash \O$, $\O^c$ is 
open. Let $\H_i$ be the maximal
interval (as subsets of the circle) 
containing $H_i$ which is a subset of $\O^c$ 
and $\H = \cup_{i=1}^r \H_i$. It is possible that several $H_i$ amalgamate
into one $\H_j$, thus $r \le p.$

Since our maps are not invertible we must modify the definition of 
topologically transitive component.  We require the the forward orbit
is dense, and that the set is forward invariant.  We gather here some
facts about topologically transitive components of continuous maps
which we will use.

\begin{proposition}
1) If a topologically transitive component contains an isolated
point, then this point is periodic and the component coincides with
this periodic orbit.

2) Topologically transitive components are finite or uncountable.

3) Any dense orbit is recurrent
\end{proposition}

\begin{proof}\label{trans}
1) Let $X$ denote the topologically transitive component.
If $z$ is isolated and the orbit of $x$ is dense then $f^ix=z$ for some
$i \ge 0$. The point $x$ is also isolated, for if not then by 
continuity $z$ is not be isolated.
If the orbit of $x$ is not periodic then since $x$ is
isolated we have $fX = X \setminus \{x\}$ and $X$ is not forward 
invariant.  This 
contradiction implies that $x$ is periodic.  Since $x$'s orbit is dense
it must coincide with the topologically transitive component.

2) A component which is not finite can not contain any isolated points.
It is a simple exercise to show that a
closed set without isolated points can not be countable.

3) If $X$ is finite then this is clear.  Suppose $X$ is uncountable and
the forward orbit of $x$ is dense in $X$. Since $x$ is not isolated
the forward orbit of $x$ must come arbitrarily close to $x$ to be dense.
\end{proof}

Let $\{X_i: X_i \subset M\}$ be the topologically transitive components of 
$f|_{\Omega}$
and $\{Y_i: Y_i \subset \Sigma\}$ the topologically transitive components of
$\sigma|_{\Sigma}$.
Let $\mathcal{I}$ be the index set of the topologically 
transitive components of $f|_{\O}$.

\begin{theorem}
Every interval exclusion system has a finitely many
topologically transitive
components, in particular the number is bounded by  the equation:
$$\sharp \{i \in \mathcal{I} : X_i \text{ at most countable}\} +
 2 \sharp \{i \in \mathcal{I} : X_i \text{ uncountable}\} 
 \le 2r + \sharp \partial \P.$$
\label{lemma}
\end{theorem}
\begin{proof}
For all $z \in \O^c$ let $G(z)$ be the maximal 
interval containing
$z$ which is a subset of $\O^c$, we refer to $G(z)$ as a gap.  Let
$n(z)$ be the smallest positive integer such that $f^{n(z)}z \in \H.$
Let $k=k(z)$ be the index such that $f^{n(z)} \in \H_{k}.$ Let $J(z)$ be 
the set of $x \in G(z)$ such that
the orbit of $x$ falls into the same component $\H_{k}$ of the
maximal hole at the same time, i.e.~$n(x) = n(z)$ and $k(x) = k(z)$. 
We claim that the continuity of $f$ implies that
$J(z) = G(z).$ Indeed, if this is not the case, then 
$\partial J(z) \subset f^{-n}\partial \H_{k}.$  Since $\partial \H \subset \O$
this contradicts the definition of the gap $G(z)$.

Until further notice 
we assume that $\S$ is a SFT.  Each topologically transitive
component $X_i$ is closed, thus we can define
$a_i := \min(X_i)$ and $b_i := \max(X_i)$.
Since $\S$ is a SFT it has only a finite number of topologically transitive
components $Y_i$ and these components are disjoint.  Until
further notice we also suppose that
all the points $a_i$ and $b_i$ have unique coding (i.e.~the orbits of 
$a_i$ and $b_i$ do not intersect $\partial \P$).  With this additional
assumption the disjointness of the $Y_i$ implies that $a_i$ and $b_i$
do not belong to any transitive component other than $X_i$.
This implies that there are gaps
$G_{a_i}$ ($G_{b_i}$) on the left (right) side of
$a_i$ ($b_i$). Fix $x \in G_{a_i}$. As we saw above all $z$ between $x$ and
$a_i$ fall into the same hole at the same time $n(z)$.
By the continuity of $f$ this implies that $f^{n(z)}a_i \in \partial \H$.
A similar statement holds for $G_{b_i}$.

Consider the ordering on $\S$ which is compatible with the ordering
on $\O$.  This ordering can always be defined in an inductive manner
by simply considering the relative order of the elements of $\P^{(n)} :=
\vee_{i=0}^n f^i \P$.
(if $f$ is locally order preserving then this is simply the lexicographical 
order on $\S$).
Fix an $i \in \mathcal{I}.$
Consider the symbolic coding of $a_i$ and $b_i$. 
Call these codings $s = (s_j)_{j \in \mathcal{N}}$ and 
$t = (t_j)_{j \in \mathcal{N}}$ (here the dependence on $i$
is suppressed since $i$ is fixed).  We claim that if $X_i$
is uncountable then $s$ is not a preimage of $t$ and vice versa $t$ is
not a preimage of $s$.  To see this fix a higher block coding which
defines    Markov transition graph of $\S$.  The fact that $a_i$ is defined 
via a minimum implies that if there are several followers of the symbol
$s_j$ in the Markov graph, 
then $s_{j+1}$ is minimal follower in the sense that 
in the ordering it is smaller than all other followers.
In a similar fashion the sequence $t$ is maximal.  

If $s$ is a preimage
of $t$ or vice versa, then $s$ and $t$ are eventually maximal and
minimal at the same time.  This means that from the point on that
they agree there the maximal follower of $s_j$ is also the minimal
follower of $s_j$, so there is only one follower of $s_j$. Thus
$s$ and $t$ are eventually periodic and the $Y_i$ is simply this
periodic orbit.  In particular $Y_i$ is finite, finishing the proof of the
claim.

If $X_i$ is uncountable, then we have just shown the disjointness
of the codes of $a_i$ and $b_i$. Since we are still assuming that
the points $a_i$ and $b_i$ 
have unique coding this implies the disjointness of their $f$--orbits.
Thus at least two point of $X_i$ lie on $\partial \hat{H}$.

On the other hand if $X_i$ is at most countable, then by proposition
\ref{trans} it is finite and consists of a single periodic orbit.  
The points $a_i$ and $b_i$ lie on this orbit.

Under the above assumptions we have shown that
for every uncountable $X_i$ there are at least two points of
$X_i$ on the border of $\H$ and for every at most
countable $X_i$ there is at least one such point.
Thus, since the endpoints can only belong to a single $X_i$ we have shown
 $$\sharp \{i \in \mathcal{I} : X_i \text{ at most countable}\} +
 2 \sharp \{i \in \mathcal{I} : X_i \text{ uncountable}\} 
 \le 2r.$$

Since $\pi$ is at most 2 to 1, if we drop the assumption that the 
$a_i$ and $b_i$ have unique coding
then (at most) two symbolically disjoint topologically transitive 
components can share such a point
when projected to the interval $I$.  This can happen at most
$\sharp \partial P$ times yielding the formula in the statement of the theorem
in the case that $\S$ is a SFT.

Finally if $\S$ is not a SFT we approximate $H$ by Markov
holes $H^{(n)}$ in the same way as in theorem \ref{prev}.
Arguing similarly to the proof of theorem \ref{prev} we have
$\tilde{\S} := \overline{\cup_j \S_j} \subset \S$ and the difference 
consists of BSCs.
Note that the number of transitive component of $\tilde{\S}$ 
is less than or equal to 
the limsup of the number of components of the approximating sequence.
It remains to bound the number of BSCs.

Let $u \le 2r$ be the number of points in $\partial \hat{H}$ 
which are boundary 
points of non BSCs. By the definition of a BSC such a point can not belong
to a BSC. 
We claim that the number of BSCs is at most $2r-u$.
Let $\mathcal{C} := \{c_j\}$ be the set of points in $\partial H$ 
which belong to a BSC. By definition of BSCs 
$\mathcal{C}$ must be disjoint from $\pi(\tilde{\S})$ or else they
would be approximable by periodic points. Also
by definition every $c_j$ must belong to $\partial \hat{H}$.
If two holes $H_i$ and $H_j$ amalgamate to a single $\hat{H}_k$ 
then two points in $\partial H$ belong to $\hat{H}_k$ and
(i.e.~their orbit falls into $H$) and thus can not belong to $\mathcal{C}$.
Thus the cardinality of $\mathcal{C}$ is at most $2r-u$. 
Clearly any BSC $X_i$
contains at least one of the points $c_j$.  We will associate with
each BSC $X_i$ a subset $I(X_i)$ of $\partial H \subset X_i$ of 
``insertable boundary
points''. We claim that two distinct BSCs must have nonintersecting 
sets $I$.  Once this claim has been established it immediately follows that
the number of BSCs is at most $2r-u$ which will complete the proof of
the theorem.

We turn to the proof of the claim. Fix $X = X_i$ a BSC.
Let $\hc_j = \pi^{-1} c_j \in \S$.  Note that $\hc_j$ is uniquely
defined by the convention we made in section 2.  
Fix $x \in X$ such that the orbit of $x$ is dense in $X$.
Consider any $\hx \in \S$ such that $\pi(\hx)= x$.
For each $k \ge 0$ we define $n_k(\hx)$ in the following way:
$n_k(\hx)$ is the longest initial block of one of the $\hc_j$ 
which agrees with the initial block of the same length of $\sigma^k\hx$.
We call $n_k(\hx)$ the flag of $x$ at time $k$.
We remark that $n_k(\hx) < \infty$ if $f^kx \not \in \partial H$. 
By proposition \ref{trans} either $x$ is an isolated periodic point which 
hits $\partial H$ or it visits $\partial H$ only a finite number of times.
In the second case since by the proposition $x$ is recurrent, we can assume
by replacing $x$ by $f^ix$ for sufficiently large $i$ that the orbit of $x$
does not visit $\partial H$ at all.

We remind that $x$ is fixed and thus we drop the $x$ dependence
of the notations.
We say that a flag $n_k$ is covered by another flag if
there is a $k' < k$ such that $n_{k'} + k' \ge n_k + k$
(see figure 1).  If a flag is not covered by any other flag we
say it is uncovered.  The motivation for this terminology is
the following fact:  if the flag of $x$ is uncovered at time $k$ then we can
concatenate $\hx$'s initial segment of length $k$ with 
any orbit in $\S$ whose initial segment agrees with 
the next $n_k$ entries of $\hx$ and the concatenated orbit belongs
to $\S^*$ (it may be wandering and thus not belong to $\S$).

We claim that there are infinitely many 
indices $t_l$ such that the flags at these times are
uncovered. To see this this consider $n_1 + 1$.
Clearly there must be a $k >1 $ such that $n_k + k > n_1 + 1$.
Let $t_1$ be the smallest such $k$, this flag must be uncovered.
Arguing inductively the flags at times 
$t_{l+1} := \min\{k > t_l: n_k + k > n_{t_l} + t_l\}$ are uncovered.

For each $k \ge 0$ we define the color(s) of $x$ at time $k$
to be the set of $c_j$
such that the initial blocks of length $n_k(\hx)$ of $\hc_j$ 
and $\sigma^k\hx$ coincide.  We assume that we start with a fine
enough Markov partition that each $c_j$ belongs to a different
element of the time zero partition.  This implies that 
the color of $x$ is unique for each $k$.
Let $I(X) := I(X,x)$ be the set of colors which occur infinitely often
at uncovered times. Clearly this set is nonempty.

We are now ready to prove the above claim.
Consider $I(X,x)$ and $I(Z,z)$ which have nonempty intersection, 
and let $c = c_j$ belong to this intersection.
We will recursively construct a point $w \in \S$ such that
the orbit of $\pi(w)$ is dense in $X \cup Z$.
To construct $w$ fix a positive sequence $\e_m \to 0$. 
Consider an initial segment of $\hx$ of length $t_{l_0}$ where
$t_{l_0}$ is the smallest $t_l(x)$ such that the orbit segment 
$f^ix: i=0, \dots, 
t_{l_0}$ is $\e_0$ dense in $X$.
Consider the smallest $k_1$ such that $\sigma^{k_1} \hz$ and $\hc$ agree at at 
least
$n_{t_{l_0}}(\hx)$ places.  Furthermore consider am initial segment of
$\sigma^{k_1}\hz$  of length $t_{l_1}$ where
$t_{l_1}$ is the smallest $t_l(f^kz)$ such that the orbit segment 
$f^iz: i=0, \dots, 
t_{l_0}$ is $\e_1$ dense in $Z$. Consider the smallest $k_2$ such that 
$\sigma^{k_2}\hx$ and $\hc$ agree at at least
$n_{t_{l_1}}(\sigma^{k_1}\hz)$ places. 
These definitions guarantee that the sequence
$\hx(t_{l_0}) \star \sigma^{k_1}\hz(t_{l_1}) \star \sigma^{k_2}\hx$
belongs to $\S^*$.
Here the symbol $\star$ mean concatenation
and $\hx(n) = \hx_0,\dots,\hx_n$.
Recursively repeating this construction produces the sequence
$$w := \hx(t_{l_0}) \star \sigma^{k_1}\hz(t_{l_1}) 
\star \sigma^{k_2}\hx(t_{l_2})
\star \sigma^{k_3}\hz(t_{l_3}) \star \cdots.
$$
By construction $w \in \S^*$ and 
the orbit of $\pi(w)$ is clearly dense in $X \cup Z$.
In particular $\pi(w)$ in nonwandering and thus $w \in \S$.
This finishes the proof of the claim.
\end{proof}

\section{Genericity results}

\subsection{A criterion for SFTs}
\begin{lemma}
If for each $x \in \bH$ there is an $i$ such that
$f^ix \in \iH$ then $f|_{\O}$ is a uniformly hyperbolic
diffeomorphism  and $\S$ is a SFT.
\label{lemma1}
\end{lemma}

\begin{proof} Fix a generating Markov partition $\P$ and let
$\P^{(n)} := \vee_{i=-n}^n f^i \P$.
Let  $P^{(n)}_x := \bigcup_{\{P \in \P^{(n)}: x \in P \}} P$. 
If $x \in \bH$ and $f^i x \in \iH$ then since $\P$ is generating by
continuity there is a $n(x)$ such
that $f^i P^{(n(x))}_x \subset \iH$.
Since $\bH$ is
compact we can cover $\bH$ by a finite collection of the sets
$P^{(n(x))}_x$ to obtain a neighborhood $N$ of $\bH$ such that 
$N \cap \O = \emptyset.$  Since we used a finite collection of $P_x$, the
hole $H' := N \cup H$ consists of a finite union of element of $\P^{(N)}$ for
some sufficiently large integer $N$ and thus $\S$ is a SFT.
\end{proof}

\subsection{Results in the Hausdorff metric}

Let $s$ denote the dimension of $M$.  Consider the set 
$C$ of all holes such that $\bH$ is a continuous, i.e.~there is a
continuous map $h: {\bf S^{s-1}} \rightarrow M$ whose image is $\bH$.  For 
$H \in C$ let 
$H^{\epsilon} := \cup_{x \in \bH} B(x,\epsilon).$
For $H_1,H_2 \in C$ we define
$$d(H_1,H_2) := \inf \big \{ \epsilon > 0: H_1 \subset H_2^{\epsilon},
H_2 \subset H_1^{\epsilon} \big \}.$$

\begin{lemma}
If the set $\O_H$ is totally disconnected then for every $\e > 0$
there is a hole $H' \in C$ with $d(H',H) < \e$ and an open neighborhood
$\mathcal{U} \subset C$ of $H'$ such that for all holes $H'' \in \mathcal{U}$
the subshift $\S$ is a SFT.
\label{totallydisconnected}
\end{lemma}

\begin{proof}
Consider the set   
$\mathcal{V} := \{H' \in C: H \subset H', \  d(H,H') < \epsilon\}.$
Since $\O_H$ is totally disconnected there are
holes $H' \in \mathcal{V}$ such that $\bH' \cap \O_H = \emptyset$. 
Since $\O_H \subset \O_{H'}$ such an 
$H'$ satisfies the requirements of lemma \ref{lemma1} and
so defines a SFT.

Next we will show that an open set $\mathcal{U}$ of holes satisfy the 
requirements of lemma \ref{lemma1}. By compactness of $\iH'$ we can find 
an open set $B \subset \iH'$ with $d(B, H') > 0$
and a positive integer $N$ such that 
for all $x \in \bH$ there exists an integer $i$,  
satisfying $|i| \le N$ such that $f^ix \in B.$  
Then, just as in 
the proof of lemma \ref{lemma1}  
there is a neighborhood $U$ of $\bH'$ such that for each $x \in U$ 
for some $i$ satisfying $|i| \le N$ we have $f^ix \in B.$  
This immediately implies that we can choose a small neighborhood
$\mathcal{U} \subset C$ of $H'$ such that for any hole $H'' \in \mathcal{U}$
the boundary of $H''$ satisfies the requirements of lemma \ref{lemma1}
and thus the corresponding shift $\S''$ is a SFT.
\end{proof}

\begin{theorem} In the two dimensional case
the set of $H \in C$ for which $\S$ is a SFT is open and dense. 
\label{theorem4}
\end{theorem}

\begin{proof}
Consider an arbitrary hole $H \in C$ and the associated invariant
set $\O_H$. 
If $x$ is a generic point in the sense that it visits (in both
forward and backwards time) any cylinder set (defined by the Markov
partition) with the correct frequency, then 
$W^s(x)$ and $W^u(x)$ completely fall into $H$ and
thus are both disjoint from $\O_H$.  Both $W^s(x)$ and $W^u(x)$
are curves which are dense in $M$, thus since $M$ is two dimensional 
the complement of their union is totally 
disconnected (i.e. no two points are in the same connected component).  
Since $\O_H$ is a subset of this set  it is also
totally disconnect. Apply Lemma \ref{totallydisconnected} finishes
the proof.
\end{proof}

\subsection{Rectangle like holes}
\label{rectholes}
In this section we assume that $M$ is two dimensional.
We call a hole a {\em rectangle like hole} if $\bH$ consists of
a finite number of curves, each of which is
a piece of a stable or unstable manifold of $f$.  We assume that at
each point where two curves meet, that one is stable, the
other is unstable.  Thus the number of corners is always even,
and if we fix an orientation of the boundary we only need to
give the coordinates of every other corner point to describe
a rectangle like hole.  If we fix the number $2n$ of corners of a 
rectangle like hole $(n \ge 2)$
then we can parameterize the set of all such rectangles by an open
subset $R^{(n)}$ of $T^{2n}.$  We will consider the Lebesgue measure
on $R^{(n)}.$

\begin{theorem}
For every $n \ge 2$, the set of rectangle like holes with $2n$ corners
for which $\S$ is a SFT is of full measure and
contains an open dense subset of $R^{(n)}.$
\label{theorem1}
\end{theorem}

\begin{proof}
Consider the set $G$ of generic points in the same sense as in the proof
of theorem \ref{theorem4}.  The set $G$ is
of full Lebesgue measure.
Suppose $H$ is a rectangle like hole with the property that
every other corner point (those which are noted in the description
of $H \in R^{(n)}$) is in $G.$  Such a hole satisfies the
requirements of lemma \ref{lemma1} and thus $\O$ is a SFT.
Clearly the set of such holes is dense and of full measure.  
The proof of openness is the same as in the proof of the previous theorem.
\end{proof}

\subsection{Polyhedral holes}

In this section we
suppose that $M$ is $s$--dimen\-sional and has a flat structure
(i.e.~$M \subset {\bf R^s}$ or $M \subset {\bf T^s}$).
We consider holes which are the interior
of arbitrary polyhedron.  Fix  the number of corners $n \ge s+1$. The 
set of $n$-gons is an open subset of $M^n$ which we denote by
$P^{(n)}$.  

\begin{theorem}
For every $n \ge s+1$ the set of polyhedral holes for which
$\S$ is a SFT has positive Lebesgue measure.
\label{theorem2}
\end{theorem}

\begin{proof}
We call a polyhedral hole $H \in P^{(n)}$ large if
the Hausdorff dimension of the associated invariant set $\O$ is
less than or equal to one.  We note that if $H \subset H'$
then $\S_{H'} \subset \S_{H}.$  Thus if $H$ is a large hole then
so is $H'.$ 

Our proof is local.  Fix an open set $B$.
Consider the set $\S_B$ of points
which never fall into $B$.  Suppose $B$ is large enough that 
$a := \text{dim}(\O_B) \le 1.$ 
This implies that $\text{dim}(\text{proj}_{\theta}(\O_B)) = a \le 1$ for
almost every $\theta$
(\cite{F}, Theorem 6.9). 
Here $\text{proj}_{\theta}$ denotes the orthogonal projection from 
$M$ onto $L_{\theta}$, the line through the origin 
in the direction $\theta \in {\bf S^{s-1}}.$
The term almost every $\theta$ refers to the Lebesgue
measure on the ${\bf S^{s-1}}$.

Consider the set of $(s-1)$--dimensional polyhedra contained in
a co-dimension one  hyperplane with normal direction $\theta$. 
Such a polyhedron plays the role of a face of a polyhedral hole.  
Let $t$ be the parameter on $L_{\theta}.$  The projection
of any such face onto $L_{\theta}$ is simply a point $t \in L_{\theta}$.
The pair of parameters $\{\theta,t\}$ determine a family of faces which
lie in a common co--dimension one hyperplane.

Suppose $\theta$ is a generic direction, 
i.e.~$\text{dim}(\text{proj}_{\theta}(\O_B)) = a \le 1$. Then for a.e.~$t$,
any polyhedral face normal to $\theta$ will be disjoint from $\O_B.$      

Now consider any polyhedron $H$ which contains $B$ in it's interior
and for which all the faces are generic in the
above sense, i.e.~they do not intersect the set $\O_B$.  
This immediately implies that the boundary $\bH$ satisfies the
requirements of lemma \ref{lemma1}. 
The set of such polyhedra
is clearly open, and locally dense, and locally of full measure.
\end{proof}

In fact we can show a bit more.  Let $\hat{P}^{(n)} \subset P^{(n)}$ 
be the set of all large polyhedral holes.  If we choose $B \in
\hat{P}^{(n)}$ in the above proof then we have shown:

\begin{theorem}
For every $n \ge s+1$ the set of  large polyhedral holes for which
$\S$ is a SFT is  of full measure and contains an open dense
subset of $\hat{P}^{(n)}$.
\label{theorem3}
\end{theorem}

\noindent
Remark: We believe that the set of polyhedral holes defining a SFT
is of full measure and contains an open and dense set.  
Our strategy of proof  can not be used since the projection onto a
line $L_{\theta}$ of a set of dimension greater than one has
positive one dimensional measure for a.e.~$\theta$.

\end{document}